\documentclass{ifacconf}
\usepackage[centertags]{amsmath}
\usepackage{graphics}
\usepackage{amscd}
\usepackage{amsfonts}
\usepackage{amssymb}
\usepackage{helvet}
\usepackage{rotating}
\usepackage{epsfig}
\usepackage[latin1]{inputenc}
\usepackage{placeins}
\usepackage{pst-all}
\usepackage{relsize}
\usepackage{longtable}
\usepackage{natbib}
\usepackage{color}
\newtheorem{theorem}{Theorem}[section]

\newtheorem{proposition}[theorem]{Proposition}

\newtheorem{assumption}[theorem]{Assumption}

\newtheorem{remark}[theorem]{Remark}
\newtheorem{example}[theorem]{Example}
\newcommand{\R}{\mathbb{R}}
\newcommand{\N}{\mathbb{N}}

\newcommand{\X}{\mathbb{X}}
\newcommand{\U}{\mathbb{U}}
\newcommand{\XX}{X}
\newcommand{\UU}{U}

\newcommand{\cS}{\mathcal{S}}
\newcommand{\cK}{\mathcal{K}}
\newcommand{\cKL}{\mathcal{KL}}

\def\argmin{\mathop{\rm argmin}}

\usepackage{framed}
\newenvironment{fshaded}{%
\MakeFramed {\FrameRestore}}%
{\endMakeFramed}
\definecolor{shadecolor}{rgb}{1.,1.,1.}%
\definecolor{framecolor}{rgb}{.0,.0,.0}%
\begin{document}
%%%%%%%%%%%%%%%%%%%%%%%%%%%%%%%%%%%%%%%%%%%%%%%%%%%%%%%%%%%%%%%%%%%%%%%%%%%%%%%%

%%%%%%%%%%%%%%%%%%%%%%%%%%%%%%%%%%%%%%%%%%%%%%%%%%%%%%%%%%%%%%%%%%%%%%%%%%%%%%%%
\begin{frontmatter}

\title{Reducing the Prediction Horizon in NMPC: An Algorithm Based Approach\thanksref{footnoteinfo}}

\thanks[footnoteinfo]{This work was supported by DFG Grant Gr1569/12 within the Priority Research Program 1305 and the Leopoldina Fellowship Programme LPDS 2009-36.}

\author[First]{J.~Pannek} 
\author[Second]{K.~Worthmann} 

\address[First]{Curtin University of Technology, Perth, 6845 WA, Australia\newline (e-mail: juergen.pannek@googlemail.com)}
\address[Second]{University of Bayreuth, 95440 Bayreuth, Germany\newline (e-mail: karl.worthmann@uni-bayreuth.de)}

\begin{abstract}
In order to guarantee stability, known results for MPC without additional terminal costs or endpoint constraints often require rather large prediction horizons. Still, stable behavior of closed loop solutions can often be observed even for shorter horizons. Here, we make use of the recent observation that stability can be guaranteed for smaller prediction horizons via Lyapunov arguments if more than only the first control is implemented. Since such a procedure may be harmful in terms of robustness, we derive conditions which allow to increase the rate at which
state measurements are used for feedback while maintaining stability and desired performance specifications. Our main contribution consists in developing two algorithms based on the deduced conditions and a corresponding stability theorem which ensures asymptotic stability for the MPC closed loop for significantly shorter prediction horizons. 
\end{abstract}

\begin{keyword}
model predictive control, stability, suboptimality estimate, algorithms, sampled data system
\end{keyword}

\end{frontmatter}
%%%%%%%%%%%%%%%%%%%%%%%%%%%%%%%%%%%%%%%%%%%%%%%%%%%%%%%%%%%%%%%%%%%%%%%%%%%%%%%%

%%%%%%%%%%%%%%%%%%%%%%%%%%%%%%%%%%%%%%%%%%%%%%%%%%%%%%%%%%%%%%%%%%%%%%%%%%%%%%%%
\section{Introduction}
\label{Section:Introduction}
%%%%%%%%%%%%%%%%%%%%%%%%%%%%%%%%%%%%%%%%%%%%%%%%%%%%%%%%%%%%%%%%%%%%%%%%%%%%%%%%

Model predictive control (MPC), sometimes also termed receding horizon control (RHC), deals with the problem of approximately solving an infinite horizon optimal control problem which is, in general, computationally intractable. To this end, a solution of the optimal control problem on a truncated, and thus finite, horizon is computed. Then, the first part of the resulting control is implemented at the plant and the finite horizon is shifted forward in time which renders this method to be iteratively applicable, see, e.g., \citet{MRA2009}.

During the last decade, theory of MPC has grown rather mature for both linear and nonlinear systems as shown in \citet{AZ2000} and \citet{RM2009}.
Additionally, it is used in a variety of industrial applications, cf. \citet{QB2003} due to its capability to directly incorporate input and state constraints. 
However, since the original problem is replaced by an iteratively solved sequence of control problems -- which are posed on a finite horizon -- stability of the resulting closed loop may be lost. In order to ensure stabilty, several modifications of the finite horizon control problems such as stabilizing terminal constraints or a local Lyapunov function as an additional terminal weight have been proposed, cf. \citet{KG1988} and \citet{ChAl1998}. 

Here, we consider MPC without additional terminal contstraints or costs, which is especially attractive from a computational point of view. Stability for these schemes can be shown via a relaxed Lyapunov inequality, see \citet{G2009} for a set valued framework and \citet{GP2009} for a trajectory based approach. 
Still, stable behavior of the MPC controlled closed loop can be observed in many examples even if these theoretical results require a longer prediction horizon in order to guarantee stability. The aim of this paper consists of closing this gap between theory and experimental experience. In particular, we develop an algorithm which enables us to ensure stability for a significantly shorter prediction horizon %in comparison to the estimates from \citet{G2009} 
and, as a consequence, reduces the computational effort in each MPC iteration significantly. 
The key observation in order to extend the line of arguments in the stability proofs proposed in \citet{G2009} is that the relaxed Lyapunov inequality holds true for the MPC closed loop for smaller horizons if not only the first element but several controls of the open loop are implemented, cf. \citet{GPSW2010}. 
However, this may be harmful in terms of robustness, cf. \citet{MS2007}.
Utilizing the internal structure of consecutive MPC problems in the critical cases along the closed loop trajectory, we present conditions which allow us to close the loop more often and maintain stability. 
Indeed, our results are also applicable in order to deduce performance estimates of the MPC closed loop. 
%guaranteeing the validity of this relaxed Lyapunov inequality in each MPC iteration requires a sufficiently long prediction horizon, see \citet{GPSW2010} for estimates.
%To this end, we introduce a lower bound on the degree of suboptimality which we wish to guarantee and extend the trajectory based suboptimality results stated in \citet{GP2009} to cover longer open loop periods. 

The paper is organized as follows: In Section \ref{Section:Problem formulation}, the MPC problem under consideration is stated. 
In the ensuing section, we extend suboptimality results from \citet{GP2009} and \citet{GPSW2010} in order to give a basic MPC algorithm which checks these performance bounds.
In Section \ref{Section:Improving Results}, we develop the proposed algorithm further in order to increase the rate at which state measurements are used for feedback and state our main stability theorem.  
In the final Section \ref{Section:Conclusion} we draw some conclusions.
Instead of a separated example, we use a numerical experiment throughout the paper which shows the differences and improvements of the presented result in comparison to known results.

%%%%%%%%%%%%%%%%%%%%%%%%%%%%%%%%%%%%%%%%%%%%%%%%%%%%%%%%%%%%%%%%%%%%%%%%%%%%%%%%
\section{Problem formulation}
\label{Section:Problem formulation}
%%%%%%%%%%%%%%%%%%%%%%%%%%%%%%%%%%%%%%%%%%%%%%%%%%%%%%%%%%%%%%%%%%%%%%%%%%%%%%%%

We consider nonlinear discrete time control systems
\begin{align}
	\label{Setup:nonlinear discrete time system}
	x(n + 1) = f(x(n), u(n)), \quad x(0) = x_0
\end{align}
with $x(n) \in \X \subset \XX$ and $u(n) \in \U \subset \UU$ for $n \in \N_0$ where $\N_0$ denotes the natural numbers including zero. The state space $\XX$ and the control value space $\UU$ are arbitrary metric spaces. As a consequence, the following results are applicable to discrete time dynamics induced by a sampled -- finite or infinite dimensional -- system, see, e.g., \citet{AGW2010a} or \citet{IK2002}. Constraints may be incorporated by choosing the sets $\X$ and $\U$ appropriately. Furthermore, we denote the space of control sequences $u: \N_0 \rightarrow \U$ by $\U^{\N_0}$. \\
Our goal consists of finding a static state feedback $u = \mu(x) \in \U$ which stabilizes a given control system of type \eqref{Setup:nonlinear discrete time system} at its unique equilibrium $x^*$. In order to evaluate the quality of the obtained control sequence, we define the infinite horizon cost functional
\begin{align}
	\label{Setup:infinite cost functional}
	J_\infty (x_0, u) = \sum\limits_{n = 0}^\infty \ell(x(n), u(n))
\end{align}
with continuous stage cost $\ell: \X \times \U \rightarrow \R_0^+$ with $\ell(x^*,0) = 0$ and $\ell(x,u) > 0$ for $x \neq x^*$. Here, $\R_0^+$ denotes the nonnegative real numbers. %which is positive definite with respect to its first argument in order to characterize the desired equilibrium $x^*$. 
The optimal value function %corresponding to \eqref{Setup:infinite cost functional} 
is given by $V_\infty(x_0) = \inf_{u \in \U^{\N_0}} J_\infty(x_0, u)$ and, for the sake of convenience, it is assumed that the infimum with respect to $u \in \U^{\N_0}$ is attained. 
Based on the optimal value function $V_\infty(\cdot)$, an optimal feedback law on the infinite horizon can be defined as
\begin{align}
	\label{Setup:infinite control}
	\mu_\infty(x(n)) := \argmin_{u \in \U} \left\{ V_\infty(x(n + 1)) + \ell(x(n), u) \right\}
\end{align}
using Bellman's principle of optimality. However, since the computation of the desired control law \eqref{Setup:infinite control} requires, in general, the solution of a Hamilton--Jacobi--Bellman equation, which is very hard to solve, we use a model predictive control (MPC) approach instead. The fundamental idea of MPC consists of three steps which are repeated at every discrete time instant:
\begin{itemize}
	\item Firstly, an optimal control $(u_0^\star, u_1^\star, \ldots, u_{N-1}^\star) \in \U^{N}$ for the problem on a finite horizon $[0, N)$, i.e.,
		\begin{align}
			\label{Setup:finite cost functional}
			J_N(x_0, u) = \sum\limits_{k = 0}^{N - 1} \ell(x_u(k; x_0), u(k))
		\end{align}
		is computed given the most recent known state $x_0$ of the system \eqref{Setup:nonlinear discrete time system}. Here, $x_u(\cdot; x_0)$ corresponds to the open loop trajectory of the prediction with control $u$ and initial state $x_0$.
	\item Secondly, the first element $\mu_N(x_0) = u_0^\star$ of the obtained sequence of open loop control values is implemented at the plant.
	\item Thirdly, the entire optimal control problem considered in the first step is shifted forward in time by one discrete time instant which allows for an iterative application of this procedure.
\end{itemize}
The corresponding closed loop costs are given by
\begin{align*}
	J_\infty (x_0, \mu_N) = \sum\limits_{n = 0}^\infty \ell(x_{\mu_N}(n), \mu_N(x_{\mu_N}(n)))
\end{align*}
where $x_{\mu_N}(\cdot)$ denotes the closed loop solution. We use
\begin{align}
	\label{Setup:open loop control}
	u_N(\cdot; x_0) = \argmin_{u \in \U^N} J_N(x_0, u)
\end{align}
to abbreviate the minimizing open loop control sequence and $V_N(x_0) = \min_{u \in \U^{N}} J_N(x_0, u)$ for the corresponding optimal value function. Furthermore, we say that a continuous function $\rho: \R_{\geq 0} \rightarrow \R_{\geq 0}$ is of class $\cK_\infty$ if it satisfies $\rho(0) = 0$, is strictly increasing and unbounded, and a continuous function $\beta: \R_{\geq 0} \times \R_{\geq 0} \rightarrow \R_{\geq 0}$ is of class $\cKL$ if it is strictly decreasing in its second argument with $\lim_{t \rightarrow \infty} \beta(r, t) = 0$ for each $r > 0$ and satisfies $\beta(\cdot, t) \in \cK_\infty$ for each $t \geq 0$. 

\begin{figure}[!ht]
	\includegraphics[width=8.5cm]{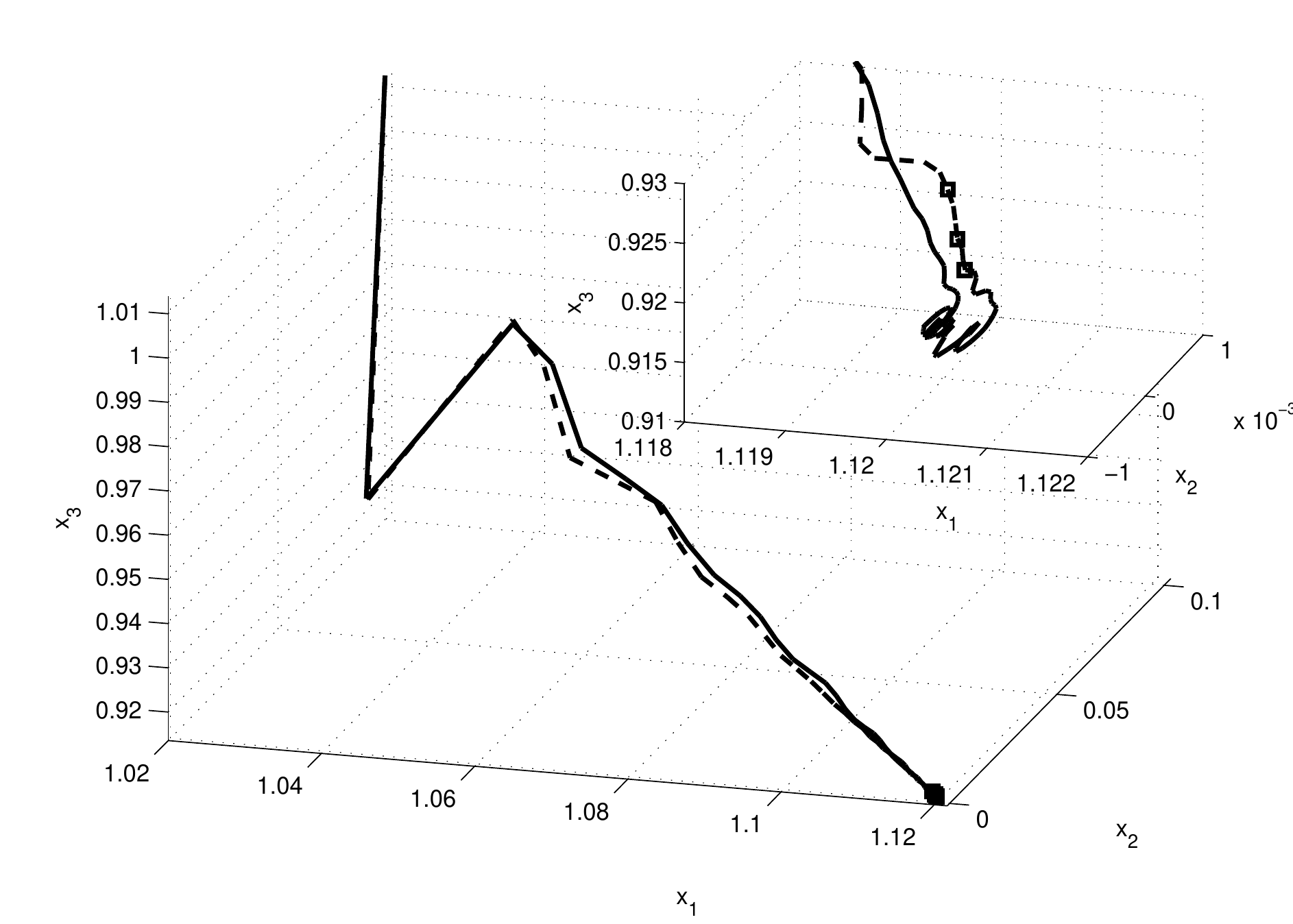}
	\caption{Comparison of closed loop solutions for $N = 30$ (solid) and $N = 19$ (dashed) where $(\square)$ marks a violation of the suboptimality bound $\overline{\alpha} = 0.1$}
	\label{Figure:comparison}
\end{figure}

Note that, different to the infinite horizon problem \eqref{Setup:infinite cost functional}, the first step of the MPC problem consists of minimizing the truncated cost functional \eqref{Setup:finite cost functional} over a finite horizon. While commonly endpoint constraints or a Lyapunov function type endpoint weight are used to ensure stability of the closed loop, see, e.g., \citet{KG1988}, \citet{ChAl1998}, \citet{JH2005} and \citet{Graichen2010}, we consider the plain MPC version without these modifications. According to \citet{GP2009}, stable behavior of the closed loop trajectory can be guaranteed using relaxed dynamic programming.
\begin{proposition}\label{Setup:proposition:aposteriori}
	(i) Consider the feedback law $\mu_N: \X \rightarrow \U$  and the closed loop trajectory $x_{\mu_N}(\cdot)$ with initial value $x(0) = x_0 \in \X$.  If
	\begin{align}
		\label{Setup:proposition:aposteriori:eq1}
		V_N(x_{\mu_N}(n)) \geq V_N(x_{\mu_N}(n + 1)) + \alpha \ell(x_{\mu_N}(n), \mu_N(x(n))) 
	\end{align}
	holds for some $\alpha \in (0, 1]$ and all $n \in \N_0$, then 
	\begin{align}
		\label{Setup:proposition:aposteriori:eq2}
		\alpha V_{\infty}(x_0) \leq \alpha J_{\infty}(x_0, \mu_N) \leq V_N(x_0) \leq V_\infty(x_0) 
	\end{align}
	holds. \\
	(ii) If, in addition, there exist $\alpha_1,\alpha_2,\alpha_3 \in \cK_\infty$ such that $\alpha_1( \| x \| ) \leq V_N(x) \leq \alpha_2( \| x \| )$ and $\ell(x,u) \geq \alpha_3( \| x \| )$ hold for all $x = x_{\mu_N}(n) \in \X$, $n \in \N_0$, then there exists $\beta \in \cKL$ which only depends on $\alpha_1, \alpha_2, \alpha_3$ and $\alpha$ such that the inequality $\| x_{\mu_N}(n) \|  \leq \beta( \| x_0 \|, n)$ holds for all $n \in \N_0$, i.e., $x_{\mu_N}(\cdot)$ behaves like a trajectory of an asymptotically stable system.
\end{proposition}
The key assumption in Proposition \ref{Setup:proposition:aposteriori} is the relaxed Lyapunov--inequality \eqref{Setup:proposition:aposteriori:eq1} in which $\alpha$ can be interpreted as a lower bound for the rate of convergence. From the literature, it is well--known that this condition is satisfied for sufficiently long horizons $N$, cf. \citet{JH2005}, \citet{GMTT2005} or \citet{AB1995}, and that a suitable $N$ may be computed via methods described in \citet{P2009b}. To investigate whether this condition is restrictive, we consider the following example.
\begin{example}\label{Setup:Example}
	Consider the synchronous generator model
	\begin{align*}
		\dot{x}_1(t) & = x_2(t) \\
		\dot{x}_2(t) & = - b_1 x_3(t) \sin x_1(t) - b_2 x_2(t) + P \\
		\dot{x}_3(t) & = b_3 \cos x_1(t) - b_4 x_3(t) + E + u(t)
	\end{align*}
	with parameters $b_1 = 34.29$, $b_2 = 0.0$, $b_3 = 0.149$, $b_4 = 0.3341$, $P = 28.22$, and $E = 0.2405$ from \citet{GORBS2003} and the initial value $x_{0} =  (1.02, 0.1, 1.014)$.
\end{example}
For Example \ref{Setup:Example} our aim is to steer the system to its equilibrium $x^{*} \approx (1.12, 0.0, 0.914)$. To this end, we use MPC to generate results for sampling time $T = 0.1$ and running cost
\begin{align*}
	\ell(x,u) = \int_0^T \| \varphi(t; x, \tilde{u}) - x^{*} \|^2 + \lambda \| \tilde{u}(t) \|^2 \, dt
\end{align*}
with piecewise constant control $\tilde{u}(t) = u$ for $t \in [0,T)$ and $\lambda = 10^{-6}$ where $\varphi(\cdot; x, \tilde{u})$ denotes the solution operator of the differential equation with initial value $x$ and control $\tilde{u}$. For this setting, we obtain stability of the closed loop via inequality \eqref{Setup:proposition:aposteriori:eq1} for $N \geq 30$, cf. Figure \ref{Figure:comparison}. However, Figure \ref{Figure:comparison} shows stable behavior of the closed loop for the smaller prediction horizon $N = 19$ without significant deviations -- despite the fact that the relaxed Lyapunov--inequality from Proposition \ref{Setup:proposition:aposteriori} is violated at three instances.
Our numerical experiment shows that, even if stability and the desired performance estimate \eqref{Setup:proposition:aposteriori:eq2} cannot be guaranteed via Proposition \ref{Setup:proposition:aposteriori}, stable and satisfactory behavior of the closed loop may still be observed. Hence, since the computational effort grows rapidly with respect to the prediction horizon, our goal consists of developing an algorithm which allows for significantly reducing the prediction horizon. To this end, we adapt the thereotical condition given in \eqref{Setup:proposition:aposteriori:eq1}. Moreover, the proposed algorithms check -- at runtime -- whether these conditions are satisfied or not.

Note that, since we do not impose additional stabilizing terminal constraints, feasibility of the NMPC scheme is an issue that cannot be neglected. In particular, without these stabilizing constraints the closed loop trajectory might run into a dead end. To exclude such a scenario, we assume the following viability condition to hold.
\begin{assumption}
	For each $x \in \X$ there exists a control $u \in \U$ such that $f(x, u) \in \X$ holds.
\end{assumption}

%%%%%%%%%%%%%%%%%%%%%%%%%%%%%%%%%%%%%%%%%%%%%%%%%%%%%%%%%%%%%%%%%%%%%%%%%%%%%%%%
\section{Theoretical Background}
\label{Section:Theoretical Background}
%%%%%%%%%%%%%%%%%%%%%%%%%%%%%%%%%%%%%%%%%%%%%%%%%%%%%%%%%%%%%%%%%%%%%%%%%%%%%%%%

In Section \ref{Section:Problem formulation} we observed a mismatch between the relaxed Lyapunov inequality \eqref{Setup:proposition:aposteriori:eq1} and the results from Example \ref{Setup:Example} which indicates that the closed loop exhibits a stable and -- in terms of performance -- satisfactory behavior for significantly smaller prediction horizons $N$. In \citet{GPSW2010} it has been shown that the suboptimality estimate $\alpha$ is increasing (up to a certain point) if more than one element of the computed sequence of control values is applied. Hence, we vary the number of open loop control values to be implemented during runtime at each MPC iteration, i.e., the system may stay in open loop for more than one sampling period, in order to deal with this issue. Doing so, however, may lead to severe problems, e.g., in terms of robustness, cf. %\citet{B1995} or 
\citet{MS2007}. Hence, the central aim of this work is the development of an algorithm which checks -- at runtime -- whether it is necessary to remain in open loop in order to guarantee the desired stability behavior of the closed loop by suboptimality arguments, or if the loop can be closed without loosing these characteristics.

\begin{figure}[!ht]
	\begin{center}
		\includegraphics[width=4.5cm]{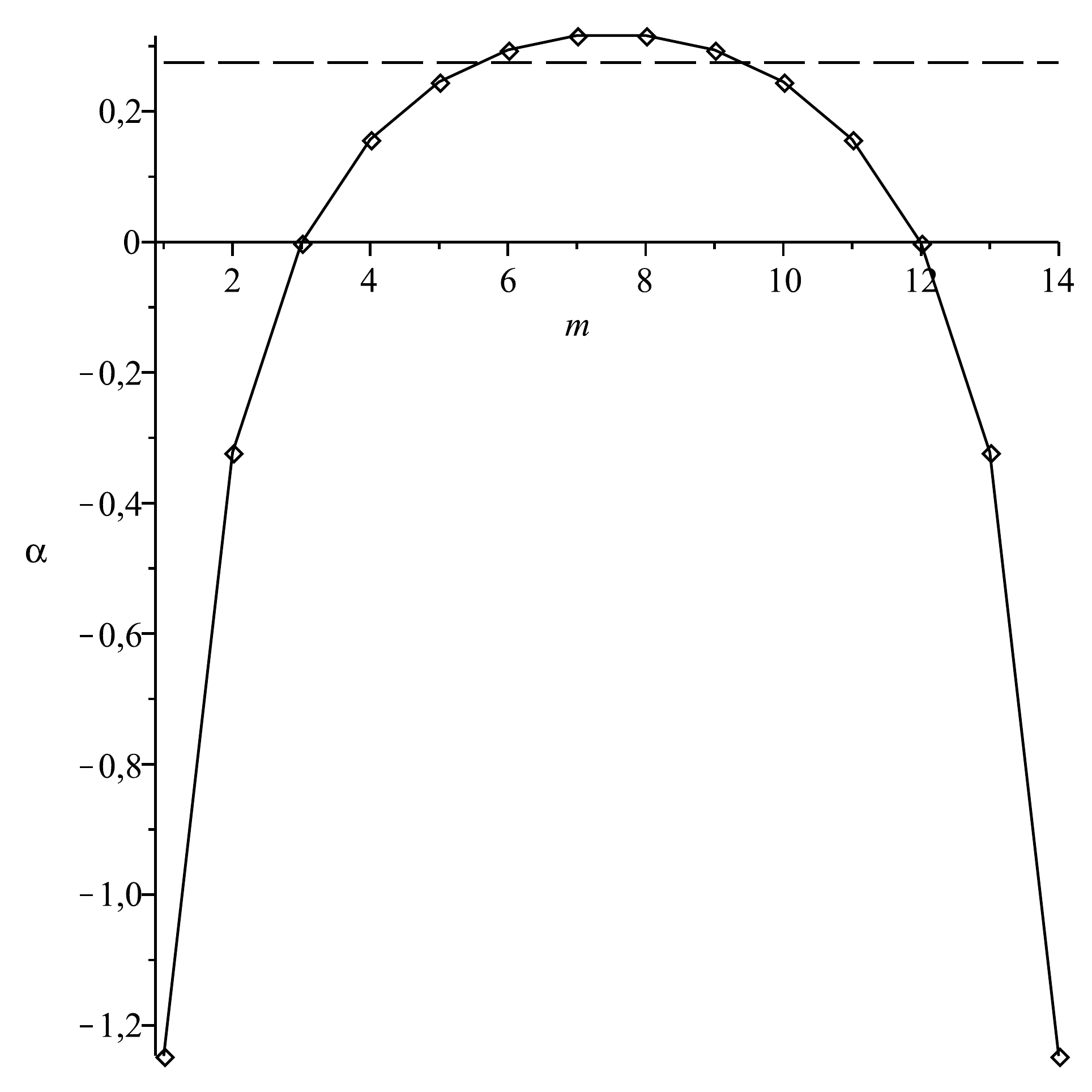}
	\end{center}
	\caption{Visualization of the $\alpha_{15,\cdot}$-values with respect to the number of controls to be implemented $m=1,\ldots,\lfloor N/2 \rfloor\}$ for prediction horizon $N=15$. The dashed line indicates the reference value $\overline{\alpha} = 0.275$.}
	\label{Theoretical Background:Figure:control horizon}
\end{figure}

\begin{figure}[!thb]
	\begin{center}
		\includegraphics[width=4.5cm]{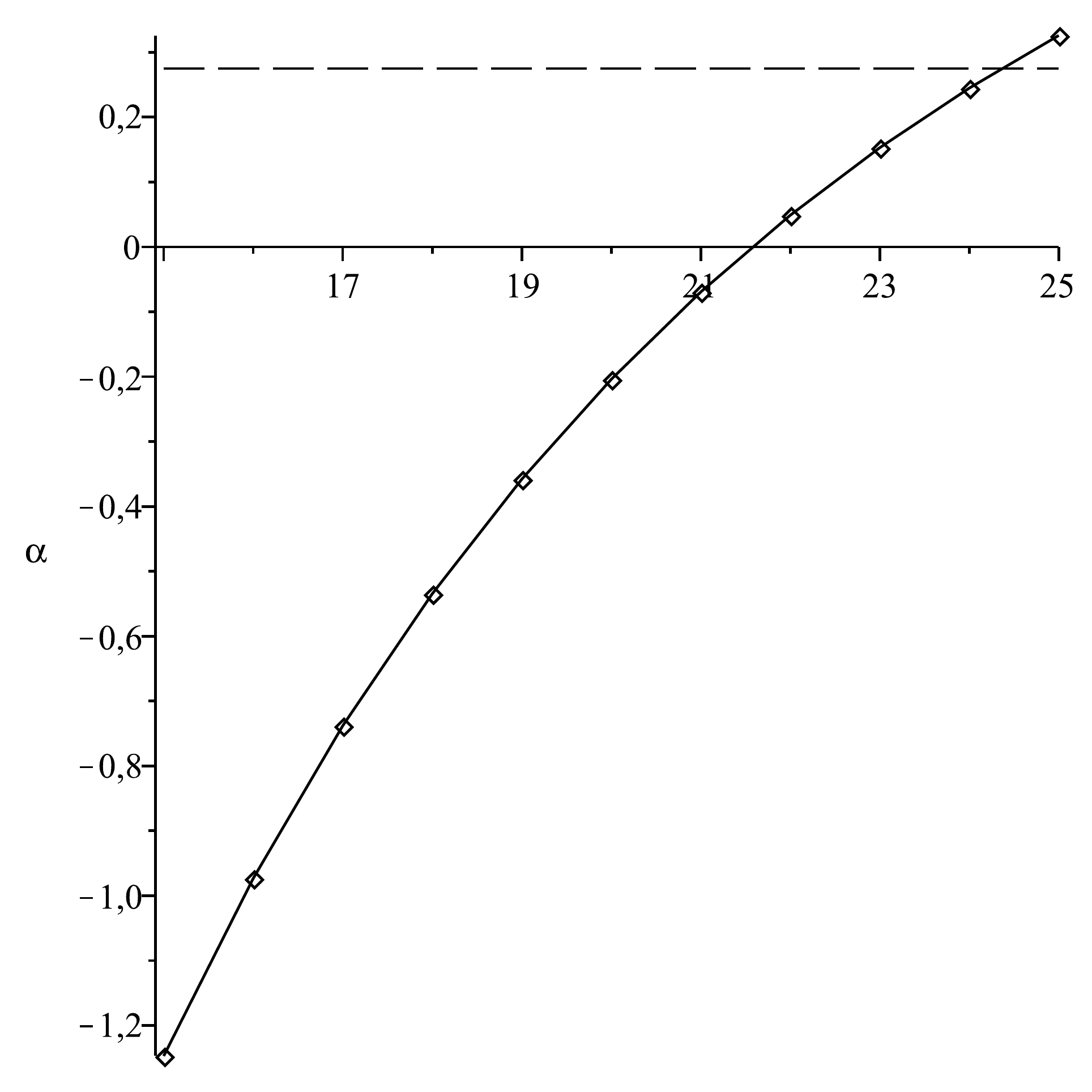}
	\end{center}
	\caption{Illustration of the $\alpha_{\cdot,1}$-values with respect to the prediction horizon $N \in \{15, \ldots, 25\}$ for $m = 1$, i.e., classical MPC. The dashed line indicates the reference value $\overline{\alpha} = 0.275$.}
	\label{Theoretical Background:Figure:prediction horizon}
\end{figure}

To this end, we introduce the list $\cS = (s(0),s(1),\ldots) \subseteq \N_0$, which we assume to be in ascending order, in order to indicate time instances at which the control sequence is updated. Moreover, we denote the closed loop solution at time instant $s(n)$ by $x_n = x_{\mu_N}(s(n))$ and define $m_n := s(n + 1) - s(n)$, i.e., the time between two MPC updates. Hence, 
\begin{equation*}
	x_{\mu_N}(s(n)+m_n) = x_{\mu_N}(s(n+1)) = x_{n+1}
\end{equation*}
holds. This enables us -- in view of Bellman's principle of optimality -- to define the closed loop control
\begin{align}
	\label{Theoretical Background:closed loop control}
	\mu_N^\cS(\cdot; x_n) & := \argmin_{u \in \U^{m_n}} \Big\{ V_{N - m_n}(x_u(m_n; x_n)) \\
	& \qquad \qquad \quad + \sum\limits_{k = 0}^{m_n - 1} \ell(x_u(k; x_n), u(k)) \Big\}. \nonumber
\end{align}
In \citet{GPSW2010}, a suboptimality degree $\alpha_{N, m_n}$ relative to the horizon length $N$ and the number of controls to be implemented $m_n$ has been introduced in order to measure the tradeoff between the infinite horizon cost induced by the MPC feedback law $\mu_N^\cS(\cdot; \cdot)$, i.e.
\begin{align}
	\label{Theoretical Background:eq:value function mu_N}
	V_\infty^{\mu_N^\cS} (x_0) := \sum\limits_{n = 0}^\infty\sum\limits_{k = 0}^{m_n - 1} \ell\left( x_{\mu_N^\cS}(k; x_n), \mu_N^\cS(k; x_n) \right),
\end{align}
and the infinite horizon optimal value function $V_\infty(\cdot)$. In particular, it has been shown that given a controllability condition, i.e., for each $x_0$ there exists a control $u_{x_0}$ such that
\begin{align}
	\label{Theoretical Background:exponential controllability assumption}
	\ell(x_{u_{x_0}}(n; x_0),u_{x_0}(n; x_0)) \leq C \sigma^n \min_{u \in \mathbb{U}} \ell(x_0,u)
\end{align}
holds with $C \geq 1$, $\sigma \in (0,1)$, then there exists a prediction horizon $N \in \mathbb{N}_{\geq 2}$ and $\overline{m} \in \{1,\ldots,\lfloor N/2 \rfloor\}$ such that $\alpha_{N,\overline{m}} \geq \overline{\alpha}$ for an arbitrarily specified $\overline{\alpha} \in (0,1)$.% Indeed, since $\alpha_{N, m}$ converges to one for $N \rightarrow \infty$ this is always possible, cf. \citet{GPSW2010}.
\begin{example}\label{Theoretical Background:Example}
	Suppose \eqref{Theoretical Background:exponential controllability assumption} to hold with $C = 4$ and $\sigma = 0.6$ and fix $\overline{\alpha} = 0.275$. Then the pair $(N,\overline{m}) = (15,6)$ implies $\alpha_{N,\overline{m}} \approx 0.294 > \overline{\alpha}$, cf. Figure \ref{Theoretical Background:Figure:control horizon}, i.e., the closed--loop satisfies the performance bound $\overline{\alpha}$. However, $N=25$ is the smallest prediction horizon which yields $\alpha_{N,m} \geq \overline{\alpha}$ for $m=1$, i.e., classical MPC, cf. Figure \ref{Theoretical Background:Figure:prediction horizon}.
\end{example}
Examples \ref{Theoretical Background:Example} shows that choosing a larger value for $m$ may significantly reduce the required prediction horizon length and, consequently, the numerical effort which grows rapidly with respect to the horizon $N$. % However, using $m > 1$ implies that the system remains in open loop for a longer period of time which may be harmful in terms of robustness.

Since the shape of the curve in Figure \ref{Theoretical Background:Figure:control horizon} is not a coincidence, i.e., $\alpha_{N, \lfloor N/2 \rfloor} \geq \alpha_{N, m}$ holds for $m \in \{1, \ldots, N - 1\}$, cf. \citet{GPSW2010}, we may easily determine the smallest prediction horizon $N \in \mathbb{N}_{\geq 2}$ such that $\alpha_{N, \lfloor N/2 \rfloor} \geq \overline{\alpha}$ ensures the desired performance specification.
However, note that these results exhibit a set--valued nature. Hence, the validity of the controllability condition \eqref{Theoretical Background:exponential controllability assumption} and the use of a prediction horizon length corresponding to the mentioned formula lead to an estimate which may be conservative -- at least in parts of the state space. This motivates the development of algorithms which use the above calculated horizon length but close the control loop more often. Here, we make use of an extension of the $m$-step suboptimality estimate derived in \citet{GP2009} which is similar to \citet{GPSW2010} but can be applied in a trajectory based setting.
\begin{proposition}\label{Theoretical Background:prop:trajectory a posteriori estimate}
	Consider $\overline{\alpha} \in (0, 1]$ to be fixed. If there exists a function $V_N: \X \rightarrow \R_0^+$ satisfying
	\begin{align}
		\label{Theoretical Background:prop:trajectory a posteriori estimate:eq1}
		V_N(x_n) \geq V_N(x_{n + 1}) + \overline{\alpha} \sum\limits_{k = 0}^{m_n - 1} \ell(x_{u_N}(k; x_n), u_N(k; x_n))
	\end{align}
	with $m_n \in \{1, \ldots, N - 1\}$ for all $n \in \N_0$, then
	\begin{align}
		\label{Theoretical Background:prop:trajectory a posteriori estimate:eq2}
		%& \overline{\alpha} V_{\infty}(x_0) \leq \overline{\alpha} V_{\infty}^{\mu_N^\cS}(x_0) \leq V_N(x_0) \leq V_\infty(x_0)
		& V_{\infty}^{\mu_N^\cS}(x_0) \leq \frac 1 {\overline{\alpha}} \cdot V_\infty(x_0)
	\end{align}
	holds for all $n \in \N_0$.
\end{proposition}
\textbf{Proof:}
	Reordering \eqref{Theoretical Background:prop:trajectory a posteriori estimate:eq1}, we obtain
	\begin{align*}
		& \overline{\alpha} \sum\limits_{k = 0}^{m_n - 1} \ell(x_{u_N}(k; x_n), u_N(k; x_n)) \leq V_N(x_n) - V_N(x_{n + 1}).
	\end{align*}
	Summing over $n \in \N_0$ yields
	\begin{align*}
		& \overline{\alpha} \sum\limits_{i = 0}^{n} \sum\limits_{k = 0}^{m_i - 1} \ell(x_{u_N}(k; x_i), u_N(k; x_i)) \\
		& \leq V_N(x_0) - V_N(x_{n + 1}) \leq V_N(x_0)
	\end{align*}
	and hence taking $n$ to infinity implies the assertion.\qed

\begin{remark}\label{remark1}
	Note that the assumptions of Proposition \ref{Theoretical Background:prop:trajectory a posteriori estimate} indeed imply the estimate
 	\begin{align}
		\label{Theoretical Background:prop:trajectory a posteriori estimate:eq3}
% 		& \overline{\alpha} V_{\infty}(x_n) \leq \overline{\alpha} V_{\infty}^{\mu_N^\cS}(x_n) \leq V_N(x_n) \leq V_\infty(x_n)
		V_{\infty}(x_n) \leq V_{\infty}^{\mu_N^\cS}(x_n) \leq \frac{1}{\overline{\alpha}} V_\infty(x_n)
	\end{align}
	for all $n \in \N_0$ which can be proven analogously.
\end{remark}
Note that $\cS$ is built up during runtime of the algorithm and not known in advance. Hence, $\cS$ is always ordered. A corresponding implementation which aims at guaranteeing a fixed lower bound of the degree of suboptimality $\overline{\alpha}$ takes the following form:
\begin{fshaded}
	Given state $x:=x_0$, list $\cS = (0)$, $N \in \mathbb{N}_{\geq 2}$, and $\overline{\alpha} \in (0,1)$
	\begin{itemize}
		\item[(1)] Set $j = 0$, compute $u_N(\cdot; x)$ and $V_N(x)$. Do
		\begin{itemize}
			\item[(a)] Set $j = j + 1$, compute $V_N(x_{u_N}(j; x))$
			\item[(b)] Compute maximal $\alpha$ to satisfy \eqref{Theoretical Background:prop:trajectory a posteriori estimate:eq1}
			\item[(c)] If $\alpha \geq \overline{\alpha}$: Set $m_n = j$ and goto 2
			\item[(d)] If $j = N - 1$: Print warning ``Solution may diverge'', set $m_n = 1$ and goto 2
		\end{itemize}
		while $\alpha < \overline{\alpha}$
		\item[(2)] For $j = 1, \ldots, m_n$ do
		\begin{itemize}
			\item[] Implement $\mu_N^\cS(j-1; x) := u_N(j-1; x)$
		\end{itemize}
		\item[(3)] Set $\cS := (\cS, \text{back}(\cS) + m_n)$, $x := x_{\mu_N^{\cS}}(m_n; x)$, goto 1
	\end{itemize}
	\vspace*{-0.1cm}
\end{fshaded}
\begin{remark}
	Here, we adopted the programming notation \textit{back} which allows for fast access to the last element of a list.
\end{remark}
\begin{remark}
	If \eqref{Theoretical Background:prop:trajectory a posteriori estimate:eq1} is not satisfied for $j \leq N - 1$, a warning will be printed out because the performance bound $\overline{\alpha}$ cannot be guaranteed. In order to cope with this issue, there exist remedies, e.g., one may increase the prediction horizon and repeat step 1. For sufficiently large $N$, this ensures the local validity of \eqref{Theoretical Background:prop:trajectory a posteriori estimate:eq1}. Unfortunately, the proof of Proposition \ref{Theoretical Background:prop:trajectory a posteriori estimate} cannot be applied in this context due to the prolongation of the horizon. Yet, it can be replaced by estimates from \citet{P2009b} or \citet{G2010} for varying prediction horizons  to obtain a result similar to \eqref{Theoretical Background:prop:trajectory a posteriori estimate:eq2}. Alternatively, one may continue with the algorithm. If there does not occur another warning, the algorithm guarantees the desired performance for $x_n$ instead of $x_0$, i.e., from that point on, cf. Remark \ref{remark1}.
\end{remark}
In order to obtain a robust procedure, however, it is preferable that the control loop is closed as often as possible, i.e., $m_n = 1$ for all $n \in \N_0$. In the following section, we give results and corresponding algorithms which -- given certain conditions -- allow us to reduce $m_n$.

%%%%%%%%%%%%%%%%%%%%%%%%%%%%%%%%%%%%%%%%%%%%%%%%%%%%%%%%%%%%%%%%%%%%%%%%%%%%%%%%
\section{Algorithms and Stability Results}
\label{Section:Improving Results}
%%%%%%%%%%%%%%%%%%%%%%%%%%%%%%%%%%%%%%%%%%%%%%%%%%%%%%%%%%%%%%%%%%%%%%%%%%%%%%%%

In many applications a stable behavior of the closed loop is observed for $m_n = 1$ even if it cannot be guaranteed using the suboptimality estimate \eqref{Setup:proposition:aposteriori:eq1}, cf. Example \ref{Setup:Example}. Here, we present a methodology to close the control loop more often while maintaining stability. In particular, if $m_n > 1$ is required in order to ensure $\alpha \geq \overline{\alpha}$ in \eqref{Theoretical Background:prop:trajectory a posteriori estimate:eq1} and, as a consequence, in step (1) of the proposed algorithm, we present conditions which allow us to update the control for some $j \in \{ 1, \ldots, m_n - 1 \}$ to
\begin{align}
	\label{Improving Results:eq:control update}
	\hat{u}_N(k; x_n) := \begin{cases}
				u_N(k; x_n), & k \leq j - 1 \\
				u_N(k - j; x_{u_N}(j; x_n)), & k \geq j
				\end{cases}
\end{align}
insert the new updating time instant $s(n) + j$ into the list $\cS$, and still be able to guarantee the desired stability behavior of the closed loop. Summarizing, our aim consists of guaranteeing that the degree of suboptimality $\overline{\alpha}$ is maintained, i.e., \eqref{Theoretical Background:prop:trajectory a posteriori estimate:eq1} is still satisfied, for the MPC law updated at a time instant $j \leq m_n - 1$ according to \eqref{Improving Results:eq:control update}.
\begin{proposition}\label{Improving Results:thm:M2one a posteriori Estimate:}
	Let $x_n \in \X$ be given and inequality \eqref{Theoretical Background:prop:trajectory a posteriori estimate:eq1} hold for $u_N(\cdot;x_n)$, $\overline{\alpha} > 0$, and $m_n \in \{2, \ldots, N - 1\}$. If the inequality
	\begin{align}
		& V_N(x_{u_N}(m_n - j; x_{u_N}(j; x_n))) - V_{N - j}(x_{u_N}(j; x_n)) \nonumber \\
		\label{Improving Results:thm:M2one a posteriori Estimate:eq1}
		& \leq ( 1 - \overline{\alpha}) \sum\limits_{k = 0}^{j - 1} \ell(x_{u_N}(k; x_n), u_N(k; x_n)) \\
		& - \overline{\alpha} \sum\limits_{k = j}^{m_n - 1} \ell(x_{u_N}(k - j; x_{u_N}(j; x_n)), u_N(k - j; x_{u_N}(j; x_n))) \nonumber
	\end{align}
	holds for some $j \in \{ 1, \ldots, m_n - 1 \}$, then the control sequence $u_N(\cdot; x_n)$ can be replaced by \eqref{Improving Results:eq:control update} and suboptimality degree $\overline{\alpha}$ is locally maintained.
\end{proposition}
\textbf{Proof:}
	In order to show the assertion, we need to show \eqref{Theoretical Background:prop:trajectory a posteriori estimate:eq1} for the modified control sequence \eqref{Improving Results:eq:control update}. Reformulating \eqref{Improving Results:thm:M2one a posteriori Estimate:eq1} by shifting the running costs associated with the unchanged control to the left hand side of \eqref{Improving Results:thm:M2one a posteriori Estimate:eq1} we obtain
	\begin{align*}
		& V_N(x_{u_N}(m_n - j; x_{u_N}(j; x_n))) - V_N(x_{u_N}(0; x_n)) \\
		& \qquad \qquad \leq - \overline{\alpha} \sum\limits_{k = 0}^{m_n - 1} \ell(x_{\hat{u}_N}(k; x_n), \hat{u}_N(k; x_n))
	\end{align*}
	which is equivalent to
	\begin{eqnarray}
		\label{Improving Results:thm:M2one a posteriori Estimate:proof:eq1}
		V_N(x_n) & \geq & V_N(x_{\hat{u}_N}(m_n; x_n)) \\
		& & + \overline{\alpha} \sum\limits_{k = 0}^{m_n - 1} \ell(x_{\hat{u}_N}(k; x_n), \hat{u}_N(k; x_n)), \nonumber
	\end{eqnarray}
	i.e., the relaxed Lyapunov inequality \eqref{Theoretical Background:prop:trajectory a posteriori estimate:eq1} for the updated control $\hat{u}_N(\cdot; x_n)$.\qed
\begin{example}\label{Improving Results:Example2}
	Consider the illustrative Example \ref{Setup:Example} with horizon $N = 19$ and $\overline{\alpha} = 0.1$. For this setting, we obtain $\alpha \geq \overline{\alpha}$ with $m_n = 1$ for all updating instances $s_n \in \cS$ except for those three points indicated by $(\square)$ in Figure \ref{Figure:comparison} where we have $m_n = 2$ once and $m_n = 3$ twice. Yet, for each of these points, inequality \eqref{Improving Results:thm:M2one a posteriori Estimate:eq1} holds with $j = 1$. % if we compute $u_N(\cdot; x_{u_N}(1; x_n))$ via MPC with initial value $x_{u_N}(1; x_n)$. 
	\begin{itemize}
		\item $m_n = 2$: We apply Theorem \ref{Improving Results:thm:M2one a posteriori Estimate:} in order to update the control $u_N(\cdot; x_n) := \hat{u}_N(\cdot; x)$ according to \eqref{Improving Results:eq:control update}. %Thus, we applied $\mu_N(x(s(n))$ and $\mu_N(x(s(n) + 1))$ instead of $\mu_N^{\cS}(i; x(s(n)))$, $i = 0, 1$, i.e., we obtained stability of the classical MPC closed loop since the MPC law was updated at each sampling instant. \\
		Thus, since the MPC law was updated at each sampling instant, we have obtained stability of the classical MPC closed loop.
		\item $m_n = 3$: We utilize Theorem \ref{Improving Results:thm:M2one a posteriori Estimate:} in an iterative manner: since \eqref{Improving Results:thm:M2one a posteriori Estimate:eq1} holds for $j = 1$, we proceed as in the case $m_n = 2$, update the control $u_N(\cdot; x_n) := \hat{u}_N(\cdot; x_n)$, and compute $x_{\mu_N}(2;x_n)$. Then, we compute $u_N(\cdot; x_{\mu_N}(2; x_n))$ and check whether or not \eqref{Improving Results:thm:M2one a posteriori Estimate:eq1} holds for \eqref{Improving Results:eq:control update} with $j = 2$. Note that here, $u_N(\cdot; x_n)$ in \eqref{Improving Results:eq:control update} is already an updated control sequence. Since \eqref{Improving Results:thm:M2one a posteriori Estimate:eq1} holds true, we obtain stability of the closed loop for the updated MPC control sequence. Moreover, we applied $\mu_N(x(s(n))$, $\mu_N(x(s(n) + 1))$ and $\mu_N(x(s(n) + 2))$ instead of $\mu_N^{\cS}(i; x(s(n)))$, $i = 0, 1, 2$, i.e., we have again used classical MPC.
	\end{itemize}
	Note that this iterative procedure can be extended in a similar manner for all $m_n \leq N - 1$.
\end{example}
Integrating the results from Theorem \ref{Improving Results:thm:M2one a posteriori Estimate:} into the algorithm displayed after Proposition \ref{Theoretical Background:prop:trajectory a posteriori estimate} can be done by changing Step (2).
\begin{fshaded}
	\vspace*{-0.3cm}
	\begin{itemize}
		\item[(2)] For $j = 1, \ldots, m_n$ do
		\begin{itemize}
			\item[(a)] \hspace{- 1.5mm} Implement $\mu_N^\cS(j-1; x) := u_N(j-1; x)$
			\item[(b)] \hspace{- 1.5mm} Compute $u_N(\cdot; x_{u_N}(j; x))$ and $V_N(x_{\hat{u}_N}(m_n; x))$
			\item[(c)] \hspace{- 1.5mm} If $j < m_n$ and condition \eqref{Improving Results:thm:M2one a posteriori Estimate:eq1} holds:
			\begin{itemize}
				\item[] Construct $\hat{u}_N(\cdot; x)$ according to \eqref{Improving Results:eq:control update}
				\item[] Update $u_N(\cdot; x) := \hat{u}_N(\cdot; x)$
			\end{itemize}
 		\end{itemize}
	\end{itemize}
	\vspace*{-0.1cm}
\end{fshaded}
Note that $V_{N - j}(x_{u_N}(j; x_n))$ in \eqref{Improving Results:thm:M2one a posteriori Estimate:eq1} is known in advance from $V_N(x_n)$ due to the principle of optimality. Hence, only $u_N(\cdot; x_{u_N}(j; x_n))$ and $V_N(x_{\hat{u}_n}(m_n; x_n))$ need to be computed. Unfortunately, this result has to be checked for all $j \in \{1, \ldots, m_n - 2\}$ using the comparison basis $V_N(x_n)$. As a consequence, we always have to keep the updating instant $s(n)$ in mind, cf. Example \ref{Improving Results:Example2}. The following result allows us to weaken this restriction:
\begin{proposition}\label{Improving Results:thm:M2one a posteriori Estimate2}
	Let $x_n \in \X$ be given and inequality \eqref{Theoretical Background:prop:trajectory a posteriori estimate:eq1} hold for $u_N(\cdot;x_n)$, $\overline{\alpha} > 0$, and $m_n \in \{2, \ldots, N - 1\}$. If the inequality
	\begin{align}
		& V_N(x_{u_N}(m_n - j; x_{u_N}(j; x_n))) - V_N(x_{u_N}(m_n; x_n)) \nonumber \\
		\label{Improving Results:thm:M2one a posteriori Estimate2:eq1}
		& \leq  \overline{\alpha} \sum\limits_{k = j}^{m_n - 1} \big( \ell(x_{u_N}(k; x_n), u_N(k; x_n)) \\
		& \qquad - \ell(x_{u_N}(k - j; x_{u_N}(j; x_n)), u_N(k - j; x_{u_N}(j; x_n))) \big) \nonumber
	\end{align}
	holds for some $j \in \{ 1, \ldots, m_n - 1 \}$, then the control sequence $u_N(\cdot; x_n)$ can be replaced by \eqref{Improving Results:eq:control update} and suboptimality degree $\alpha$ is locally maintained.
\end{proposition}
\textbf{Proof:}
% 	We reformulate the relaxed Lyapunov condition
% 	\begin{align*}
% 		& V_N(x_{u_N}(m_n; x_n)) - V_N(x_n) \leq \\
% 		& \qquad \leq - \overline{\alpha} \sum\limits_{k = 0}^{m_n - 1} \ell(x_{u_N}(k; x_n), u_N(k; x_n)).
% 	\end{align*}
% 	Adding \eqref{Improving Results:thm:M2one a posteriori Estimate2:eq1} to this inequality leads to
% 	\begin{align*}
% 		& V_N(x_{\hat{u}_N}(m_n; x_n)) - V_N(x_n) \\
% 		& \qquad \qquad \leq  - \overline{\alpha} \sum\limits_{k = 0}^{j - 1} \ell(x_{u_N}(k; x_n), u_N(k; x_n)) \\
% 		& \qquad \qquad \qquad - \overline{\alpha} \sum\limits_{k = j}^{m_n - 1} \ell(x_{\hat{u}_N}(k; x_n), \hat{u}_N(k; x_n)).
% 	\end{align*}
% 	Updating the control sequence according to \eqref{Improving Results:eq:control update}, we obtain \eqref{Improving Results:thm:M2one a posteriori Estimate:proof:eq1} which completes the proof.\qed
	Adding \eqref{Improving Results:thm:M2one a posteriori Estimate2:eq1} to \eqref{Theoretical Background:prop:trajectory a posteriori estimate:eq1} leads to
	\begin{align*}
		& V_N(x_{\hat{u}_N}(m_n; x_n)) - V_N(x_n) \\
		& \qquad \qquad \leq  - \overline{\alpha} \sum\limits_{k = 0}^{j - 1} \ell(x_{u_N}(k; x_n), u_N(k; x_n)) \\
		& \qquad \qquad \quad - \overline{\alpha} \sum\limits_{k = j}^{m_n - 1} \ell(x_{\hat{u}_N}(k; x_n), \hat{u}_N(k; x_n)).
	\end{align*}
	Updating the control sequence according to \eqref{Improving Results:eq:control update}, we obtain \eqref{Improving Results:thm:M2one a posteriori Estimate:proof:eq1} which completes the proof.\qed
\begin{remark}
	We also like to point out that conditions \eqref{Improving Results:thm:M2one a posteriori Estimate:eq1} and \eqref{Improving Results:thm:M2one a posteriori Estimate2:eq1} allow for a less fast decrease of energy along the closed loop, i.e., the case $V_N(x_{\hat{u}_N}(m_n; x_n)) \geq V_N(x_{u_N}(m_n; x_n))$ is not excluded in general.
\end{remark}
As outlined in Example \ref{Improving Results:Example2}, Theorem \ref{Improving Results:thm:M2one a posteriori Estimate:} (and analogously Theorem \ref{Improving Results:thm:M2one a posteriori Estimate2}) may be applied iteratively. Yet, in contrast to \eqref{Improving Results:thm:M2one a posteriori Estimate:eq1}, $V_N(x_n)$ is not required in \eqref{Improving Results:thm:M2one a posteriori Estimate2:eq1}. Hence, if an appropriate update can be found for time instant $j \in \{1, \ldots, m_n - 1\}$, the loop can be closed and, as a consequence, Theorem \ref{Improving Results:thm:M2one a posteriori Estimate2} can be applied to the new initial value $x_{\hat{u}}(j; x_n)$ with respect to the reduced number of control to be implemented $m_n - j$. To this end, the following modification of Step (2c) can be used:

\begin{fshaded}
	\vspace*{-0.3cm}
	\begin{itemize}
% 		\item[(2)] For $j = 1, \ldots, m$ do
% 		\begin{itemize}
% 			\item[(a)] \hspace{- 1.5mm} Implement $\mu_N^\cS(j-1; x) := u_N(j-1; x)$
% 			\item[(b)] \hspace{- 1.5mm} Compute $\overline{u}_N(\cdot; x_{u_N}(j; x))$ and $V_N(x_{\hat{u}_N}(m; x))$
% 			\item[(c)] \hspace{- 1.5mm} If $j < m$ and condition \eqref{Improving Results:thm:M2one a posteriori Estimate2:eq1} holds:
			\item[(2c)] \hspace{- 1.5mm} If $j < m_n$ and condition \eqref{Improving Results:thm:M2one a posteriori Estimate2:eq1} holds:
			\begin{itemize}
				\item[] Construct $\hat{u}_N(\cdot; x)$ according to \eqref{Improving Results:eq:control update}
				\item[] Update $u_N(\cdot; x) := \hat{u}_N(\cdot; x)$
				\item[] Update the values $V_N(x_{u_N}(m_n; x))$ and \\ $\ell(x_{u_N}(k; x), u_N(k; x))$, $k=j, \ldots, m_n - 1$
			\end{itemize}
% 		\end{itemize}
	\end{itemize}
	\vspace*{-0.1cm}
\end{fshaded}
% Here, we like to point out that it is possible to choose different values $\overline{\alpha} > 0$ in Steps (1) and (2) of such an algorithm. The overall guaranteed degree of suboptimality in this case is given by the minimum of the used bounds. The positive effect of such a setting is that conditions \eqref{Improving Results:thm:M2one a posteriori Estimate:eq1} and \eqref{Improving Results:thm:M2one a posteriori Estimate2:eq1} are more likely to be satisfied.
Combining the results from Proposition \ref{Theoretical Background:prop:trajectory a posteriori estimate} and Theorems \ref{Improving Results:thm:M2one a posteriori Estimate:}, \ref{Improving Results:thm:M2one a posteriori Estimate2}, we obtain our main stability result:
\begin{theorem}
	Let $x_0 \in \X$, $\overline{\alpha} > 0$ be given and apply one of the three proposed algorithms. Assume that (1c) is satisfied for some $j \in \{1,\ldots,N-1\}$ for each iterate. Then the closed loop trajectory satisfies the performance estimate \eqref{Theoretical Background:prop:trajectory a posteriori estimate:eq2} from Proposition \ref{Theoretical Background:prop:trajectory a posteriori estimate}. If, in addition, the conditions of Proposition \ref{Setup:proposition:aposteriori}(ii) hold, then $x_{\mu_N}(\cdot)$ behaves like a trajectory of an asymptotically stable system.
\end{theorem}
\textbf{Proof:}
	The algorithms construct the set $\cS$. Since (1c) is satisfied for some $j \in \{1, \ldots, N - 1\}$ the assumptions of Proposition \ref{Theoretical Background:prop:trajectory a posteriori estimate}, i.e., inequality \eqref{Theoretical Background:prop:trajectory a posteriori estimate:eq1} for $x_n = x_{\mu_N}(s(n))$, are satisfied which implies \eqref{Theoretical Background:prop:trajectory a posteriori estimate:eq2}. To conclude asymptotically stable behavior of the closed loop, standard direct Lyapunov techniques can be applied. \qed

\begin{example}\label{Improving Results:Example3}
	Again consider Example \ref{Setup:Example} with $N = 19$ and $\overline{\alpha} = 0.1$. Since \eqref{Improving Results:thm:M2one a posteriori Estimate2:eq1} is more conservative than \eqref{Improving Results:thm:M2one a posteriori Estimate:eq1}, the closed loop control and accordingly the closed loop solution $x_{\mu_N}(\cdot)$ will not coincide in general. This also holds true in our example: While the case $m_n = 2$ can be resolved similarly using the algorithm based on \eqref{Improving Results:thm:M2one a posteriori Estimate2:eq1}, no update of the control $u_N$ is performed for the two cases $m_n = 3$.  Here, this does not cause any change in the list $\cS$. However, in general, $\cS$ is likely to be different for the presented algorithms.
\end{example}
Since the proposed algorithms guarantee stability for the significantly smaller horizon $N = 19$, the computing time for Example \ref{Setup:Example} is reduced by almost $50\%$ -- despite the additional computational effort for Step (2).

\section{Conclusion}
\label{Section:Conclusion}
%%%%%%%%%%%%%%%%%%%%%%%%%%%%%%%%%%%%%%%%%%%%%%%%%%%%%%%%%%%%%%%%%%%%%%%%%%%%%%%%

Based on a detailed analysis of a relaxed Lyapunov inequality, we have deduced conditions in order to close the gap between theoretical stability results and stable behavior of the MPC closed loop. Using this additional insight, we developed two MPC algorithms which check these conditions at runtime. Although the additional computational effort is not negligible, the proposed algorithms may allow us to shorten the prediction horizon used in the optimization. This reduces, in general, the complexity of the problem significantly which may imply a reduction of the overall computational costs.

% %%%%%%%%%%%%%%%%%%%%%%%%%%%%%%%%%%%%%%%%%%%%%%%%%%%%%%%%%%%%%%%%%%%%%%%%%%%%%%%%
% \bibliography{../bibdatei}
% %%%%%%%%%%%%%%%%%%%%%%%%%%%%%%%%%%%%%%%%%%%%%%%%%%%%%%%%%%%%%%%%%%%%%%%%%%%%%%%%

\end{document}